\documentclass[12pt,a4paper,twoside]{article}
\usepackage{bcc2003,url,graphicx,color,paralist,calc,rotating,psfrag,overpic}

\bcctitle{Counting Lattice Triangulations}
\bccname{Volker Kaibel and G\"unter~M. Ziegler}
\bccaddress{Institut f\"ur Mathematik, MA 6--2\\
Technische Universit\"at Berlin\\
D--10623~Berlin\\
Germany}
\bccemail{\{kaibel,ziegler\}@math.tu-berlin.de}

\newcommand{\R}{\mathbb{R}}
\newcommand{\Z}{\mathbb{Z}}
\newcommand{\N}{\mathbb{N}}
\newcommand{\order}[1]{\mbox{O}({#1})}

\DeclareMathOperator{\littleo}{o}
\newcommand{\conv}{\mbox{\rm conv}}

\newcommand{\freg}{f^{\mathrm{reg}}}
\newcommand{\firreg}{f^{\mathrm{irreg}}}
\newcommand{\creg}{c^{\mathrm{reg}}}
\newcommand{\abov}{\prec}
\newcommand{\cplex}{\textsf{CPLEX}~\texttt{6.6.1}}
\newcommand{\topcom}{\textsf{TOPCOM}}

\newmathop{\deven}{{\rm DE}}
\newmathop{\dodd}{{\rm DO}}
\newmathop{\aut}{{\rm Aut}}
\newmathbin{\wreath}{{\rm wr}}

\begin{document}
\makebcctitle

\begin{abstract}
We discuss the problem
to count, or, more modestly, to estimate the number $f(m,n)$
of unimodular triangulations of the planar grid of size~$m\times n$.

Among other tools, we employ recursions
that allow one to compute the (huge) number of 
triangulations for small~$m$ and rather large~$n$ by
dynamic programming;
we show that this computation can be done in polynomial
time if $m$ is fixed, and present computational results
from our implementation of this approach.

We also present new upper and
lower bounds for large $m$ and~$n$, and we 
report about results obtained from a computer simulation of the random
walk that is generated by flips.
\end{abstract}

\section{Introduction}
\label{sec:intro}

An innocent little combinatorial counting problem
asks for the number of triangulations of a
finite grid of size $m\times n$.  That is, for $m,n\ge1$
we define $P_{m,n}:=\{0,1,\ldots,m\}\times\{0,1,\ldots,n\}$, ``the
grid''. Equivalently, the point configuration $P_{m,n}$ 
consists of all points of the integer lattice $\Z^2$
in the lattice rectangle $\conv(P_{m,n})=[0,m]\times[0,n]$
of area~$mn$.
Every triangulation of this rectangle point set that uses all the
points in $P_{m,n}$ has $(m+1)(n+1)=|P_{m,n}|$ vertices, 
$2mn$ facets/triangles,
and $3mn+m+n$ edges, $2(m+n)$ of them on the boundary, 
the other $3mn-m-n$ ones in the interior.
All the triangles are minimal lattice triangles of area $\frac12$
(that is, of determinant~$1$), which are referred to as
\emph{unimodular} triangles. The grid triangulations that
use all the points are thus called \emph{unimodular triangulations}.
The number of unimodular triangulations of the grid $P_{m,n}$
 will be denoted by $f(m,n)$.
 
 As an example, our first figure shows one unimodular triangulation of
 the $5\times 6$ grid:
\begin{center}
\input{a_triangulation3.pstex_t}
\end{center}
To get started, one notes that of course $f(m,n)=f(n,m)$,
one discovers with pleasure that $f(1,n)=\binom{2n}n$, and
one works out by hand with a bit of pain that $f(2,2)=64$ and $f(2,3)=852$.
Furthermore, one observes that the special triangulations that
decompose into $m$ vertical strips of width~$1$ yield the
lower bound
\begin{equation}
  \label{eq:easylow}
  f(m,n)\ \ge\ f(1,n)^m\ =\ \binom{2n}n^m\;
\end{equation}
so for larger $m$ and $n$ the numbers $f(m,n)$ get huge \emph{very soon};
for example, the bound yields $f(5,6) > 6 \cdot 10^{14}$.

It is equally interesting to study/enumerate more
general types of triangulations, such as
triangulations of finite point sets that do not necessarily use 
all the points, triangulations of general convex or non-convex
lattice polygons, triangulations of general position
point sets, triangulations of higher-dimensional point sets, etc.
None of these will appear here, but we refer interested readers to
Lee \cite{Lee} and De Loera, Rambau \& Santos \cite{DLRS}.

Lattice triangulations are basic combinatorial objects,
and they are fundamental discrete geometric structures;
so it is no surprise that they appear in various computational geometry
contexts. However, lattice triangulations have also been
studied intensively from different algebraic geometry angles. So, 
triangulations of a convex lattice polygon
\begin{compactitem}
\item[$\bullet$] provide
the data for Viro's \cite{Viro} famous construction method
of plane algebraic curves with prescribed combinatorics and topology,
related to Hilbert's sixteenth problem;
\item[$\bullet$]
appear in Gel'fand--Kapranov--Zelevinsky's \cite{GKZ} theory of
discriminants, where ``regular'' triangulations are in bijection with 
the vertices of the secondary polytope of the point configuration, and
\item[$\bullet$]
model torus-equivariant 
crepant resolution of singularities for toric three-folds,
where ``regular'' triangulations correspond to projective desingularizations;
see e.\,g., Kempf et~al.\ \cite[Chap.~3]{KKMS} and Dais \cite{Dais-3dToric}.
\end{compactitem}
The last two points pose the problem of counting or estimating the
number $\freg(m,n)$ of \emph{regular} triangulations of~$P_{m,n}$,
that is, of triangulations of $\conv(P_{m,n})$ whose triangles are the
domains of linearity for a piecewise linear convex function
(\emph{lifting function}); see also, for example, Sturmfels
\cite[Chap.~8]{Sturmfels-GroePoly}, Ziegler \cite[Lect.~5]{Zpoly}.
Motivated by the toric variety considerations, one would like to know
whether/that for large $m$ and~$n$ most triangulations are
non-regular.  There is no proof, but a lot of evidence that this
should be true (see, e.g., Table~\ref{tab:large}).

In this context, we must admit that despite the effort that has been
put into studying this question (see e.\,g.\ Hastings
\cite[Chap.~2]{Hastings}), it has not become clear what non-regular
triangulations really ``look like.'' The ``mother of all examples'' is
the whirlpool triangulation of the $3\times 3$ grid:
\begin{center}
  \begin{overpic}[scale=2]{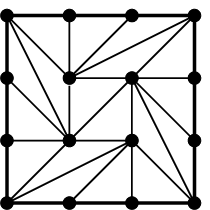}
    \put(37,29){$p_0$}
    \put(54,38){$p_1$}
    \put(54,68){$p_2$}
    \put(37,60){$p_3$}
    \put(-8,3){$q_0$}
    \put(100,3){$q_1$}
    \put(100,93){$q_2$}
    \put(-8,93){$q_3$}
  \end{overpic}
\end{center}

Indeed, this unimodular triangulation is irregular, since a lifting
function~$h$ would have to satisfy
$h(p_i)+h(q_{i-1})<h(p_{i-1})+h(q_i)$ for $i=0,\dots,3$ (all indices
taken modulo~$4$), implying the contradiction $h(q_0)+\dots+
h(q_3)<h(q_0)+\dots+ h(q_3)$.

For $m=n=3$, except for three symmetric copies
this is the only example of a non-regular triangulation: We
have $f(3,3)-\freg(3,3)=4$.
This may lead one to the conjecture that some kind of
``generalized whirlpools'' are responsible for non-regularity.
However, the pictures of ``non-regular triangulations''
that we present below do not support this intuition.
\medskip

The plan for this paper is as follows:
In Section~\ref{sec:values} we face the challenge to
\emph{\bfseries count} explicitly, trying to cope with the
``combinatorial explosion.''
For this, we present a simple dynamic programming technique by which we
get surprisingly far, and which on the specific problem of
grid triangulations outruns the much more
sophisticated general techniques such as 
Avis \& Fukuda's \cite{AvisFukuda} reverse search, 
Aichholzer's \cite{Aichholzer-path} path of a triangulation,
or the oriented matroid technique of Rambau \cite{Rambau-TOPCOM}.

In Section~\ref{sec:bounds} we \emph{\bfseries estimate} $f(m,n)$ for
large $m$ and $n$, trying to narrow the bounds for the asymptotics.
It is interesting to compare with the situation for $N=(m+1)(n+1)$
points in the plane in general position, where the currently best
available upper bound seems to be Santos \& Seidel's
\cite{SantosSeidel} estimate that there are $o(59^N)$ triangulations.
However, in our problem the $N$ points are not at all in general
position --- and the upper bounds that we have to offer are much
better: We report on a neat $O(2^{3N})$ upper bound by Anclin
\cite{Anclin}, which substantially improves on a previous $O(2^{4N})$
upper bound by Orevkov \cite{Orevkov}.  Based on explicit enumeration
results from Section~\ref{sec:values}, we get a lower bound of
$2^{2.055\,mn}$ when both $m$ and~$n$ get large; note
that~(\ref{eq:easylow}) yields already a lower bound
$2^{(1-\littleo(1))2mn}$.

Finally, in Section~\ref{sec:trias} 
we \emph{\bfseries sample} lattice triangulations for large parameters
$m$ and~$n$,
and thus try to understand what typical lattice triangulations,
as well as typical regular lattice triangulations, ``look like.''
While we have some pictures to offer, proofs seem harder to come by.
Indeed, the pictures display some long-range order;
while this may make lattice triangulations interesting as a
statistical physics model, it generates serious obstacles for any proof that
the obvious Markov chain is rapidly mixing, and thus to
application of the (by now) standard theory \cite{Behrends-MarkovChains}.


\section{Explicit Values}
\label{sec:values}

There are several methods available to generate or count all
triangulations of a finite set of points in~$\R^2$. An approach that
works for point sets in arbitrary dimensions is implemented in the
software package \topcom\ by Rambau~\cite{Rambau-TOPCOM} (see also
Pfeifle and Rambau~\cite{PfeifleRambau}). It enumerates all
triangulations in a purely combinatorial manner after the chirotope of
the point set (oriented matroid data) has been computed.

The \emph{reverse search algorithm} proposed by Avis and
Fukuda~\cite{AvisFukuda} is a rather general enumeration scheme that
can be specialized to triangulations of point configurations in $\R^2$
(see also Bespamyatnikh~\cite{Bespamyatnikh}). Since it was used to
obtain some of the results reported in Section~\ref{sec:trias}, and
because it is based upon some structural properties that are relevant
for our treatment later, we briefly describe the method here.

Let~$\mathcal{T}$ be a \emph{fine triangulation} of a point set
$S\subset\R^2$; i.e., $\mathcal{T}$ is a set of triangles, for which
the set of vertices equals~$S$, such that the union of all triangles
is the convex hull of~$S$, and any two triangles intersect in a common
(possibly empty) face. The unimodular triangulations of $P_{m,n}$ are
precisely its fine triangulations. An edge of some triangle
in~$\mathcal{T}$ is \emph{flippable} if it is contained in two
triangles of~$\mathcal{T}$ whose union is a strictly convex
quadrangle.  Replacing these two triangles by the two triangles into
which the other diagonal cuts that quadrangle (\emph{flipping} the
edge) yields another fine triangulation of~$S$.  The graph on the fine
triangulations of~$S$ defined via flipping is the 
\emph{flip graph} of~$S$.

Let us fix an arbitrary ordering of the points in~$S$, inducing via
lexicographical ordering a total order on the set of pairs of points,
and thus, again via lexicographical ordering, a total order on the
fine triangulations of~$S$, which are identified with their sets of
edges here.   There is a distinguished fine
triangulation~$\mathcal{T}_0$ of~$S$ with respect to that ordering,
namely the smallest \emph{Delaunay triangulation} (i.e., 
a triangulation characterized by the condition that for every triangle
the circumcircle contains no point from~$S$ in its interior).

Furthermore, for each fine triangulation
$\mathcal{T}\not=\mathcal{T}_0$, there is a distinguished flippable
edge
(computable in $\order{|S|}$ steps) such that, starting from any fine
triangulation of~$S$, iterated flipping of the respective
distinguished edges eventually yields $\mathcal{T}_0$.  This
algorithmically defines a spanning tree in the flip graph of~$S$,
rooted at~$\mathcal{T}_0$; in particular, the flip graph of a two
dimensional finite point set is connected.

The basic idea of the reverse search method is to traverse that
spanning tree from its root~$\mathcal{T}_0$. At each iteration one
chooses a leaf~$\mathcal{T}$ of the current partial tree, and
determines those among the triangulations adjacent to~$\mathcal{T}$ on
whose path to~$\mathcal{T}_0$ in the spanning tree~$\mathcal{T}$ lies.
Properly implemented, the reverse search algorithm generates all fine
triangulations in $\order{|S|\cdot f(S)}$ steps, where~$f(S)$ is the
number of fine triangulations of~$S$; see~\cite{AvisFukuda}.

Via the ``secondary polytope'' of a finite point set
$S\subset\R^2$~\cite{GKZ}  
one can design a variant of the reverse
search algorithm that generates \emph{all} \emph{regular}
triangulations of~$S$ in $\order{|S|\cdot F^{\text{reg}}(S)}$ steps,
where $F^{\text{reg}}(S)$ is the number of such triangulations. It is,
unclear, however, if one can also generate all \emph{regular}
\emph{fine} triangulations of~$S$ in a number of steps that is bounded
by a polynomial in the number of such triangulations.
 
If one is interested in the number of triangulations of a
two-dimensional set of points rather than in the explicit generation
of all of them, then the \emph{path of a triangulation method} due to
Aichholzer~\cite{Aichholzer-path} is more efficient than the reverse
search algorithm.

For counting the unimodular (fine) triangulations of the very special
point sets $P_{m,n}$, however, different methods are much more
efficient. These are described in the following sections.

\subsection{Narrow Strips}
\label{subsec:narrow}

\paragraph{Strips of width \boldmath$m=1$. }

For any lattice trapezoid of width~$1$, whose parallel vertical sides
have lengths $a$ and~$b$, the number of unimodular
triangulations is $g_1(a,b)=\binom{a+b}a= \binom{a+b}b$; indeed, a
bijection between these triangulations and the $a$-subsets of
$\{1,\dots,a+b\}$ is established by top-down numbering the triangles
of a triangulation~$\mathcal{T}$ by $1,\dots,a+b$, and
mapping~$\mathcal{T}$ to the $a$-set of all numbers of triangles whose
vertical edges are on the left. In particular, we have
\begin{equation}
  \label{eq:1n}
  f(1,n)\ =\ \binom{2n}{n}.  
\end{equation}

\paragraph{Strips of width \boldmath$m=2$. }

For $f(2,n)$ we have no explicit formula, and we cannot evaluate the
asymptotics precisely, but we have a ``quadratic'' recursion that can
be evaluated efficiently: For this we enumerate the triangulations
according to the highest ``width~$2$ diagonal,'' which (if it exists)
decomposes the rectangle into two width~$1$ strips, a single triangle,
and a trapezoid of width~$2$ (see the left-hand figure below).  Let
$g_2(A,B)$ denote the number of triangulations of a trapezoid of
width~$2$ with vertical edges of lengths~$A$ and~$B$ 
and horizontal base line, as in the right-hand figure below --- where
$A+B\equiv1\bmod2$ implies that the midpoint of the diagonal is not a
lattice point, and where we may assume $A<B$ by symmetry:
\[
\input{v_strip2recursion1+2.pstex_t}
\]
Thus we get 
\[
f(2,n)\ =\ \binom{2n}{n}^2\ +\ 2
\hspace{-3mm} \sum_{0\le A<B\le n\atop A+B\equiv1\bmod2}\hspace{-4mm}
g_2(A,B)
\binom{2n-\tfrac{3A+B+1}2}{n-A}
\binom{2n-\tfrac{A+3B+1}2}{n-B}.
\]
The binomial coefficients in this recursion
correspond to triangulations of width~$1$ strips. So they
could be rewritten in terms of $g_1$, as 
$g_1(n-A,n-\tfrac{A+B+1}2)$ resp.\ $g_1(n-B,n-\tfrac{A+B+1}2)$.
A similar remark applies to the binomial coefficients that
appear in the following.

For $g_2(A,B)=g_2(B,A)$ we also get a recursion by considering the
highest diagonal of width~$2$:
\begin{eqnarray*}
g_2(A,B) & = &\ 
\binom{\frac{3A+B-1}2}{A}
\binom{\frac{A+3B-1}2}{B}\\
&&+ 
\sum\limits_{0\le a\le A,\ 0\le b\le B \atop a+b\equiv1\bmod2,\ a+b<A+B}
g_2(a,b)\,
\binom{\frac{3A+B-3a-b}{2}-1}{A-a}
\binom{\frac{A+3B-a-3b}{2}-1}{B-b}.
\end{eqnarray*}
\begin{center}
\input{v_strip2recursion3+4.pstex_t}
\end{center}
Here the parameters $A,B,a,b$ may be interpreted as the $y$-coordinates
of certain lattice points. The shaded parts of the figure consist
of strips of width~$1$, whose triangulations are counted by
the binomial coefficients $g_1(\cdot,\cdot)$.

\paragraph{Strips of width \boldmath$m=3$. }

For $f(m,3)$ we have a recursion of order~$4$; it relies on the
observation that if we screen the middle strip from the 
top for diagonals of width at least~$2$, then the first diagonal
to find will be of width exactly~$2$, since any width~$3$ diagonal
is flippable, and contained in a parallelogram that is
bounded by two width~$2$ diagonals. 
The corresponding decomposition of our rectangle 
is indicated in the left figure below: 
\begin{center}
\centerline{\input{v_strip3recursion1+2.pstex_t}}
\end{center}
Therefore, we obtain
\[
f(3,n) \ =\ \binom{2n}{n}^3 + 2
\sum_{0\le A,B\le n\atop A+B\equiv1\bmod2}h(A,B,n,n)
\binom{2n-\frac{3A+B+1}{2}}{n-A}
\binom{2n-\frac{A+3B+1}{2}}{n-B},
\]
where $h(A,B,C,D)$ counts the number of triangulations
in a ``hook'' shape as given in the right drawing in the figure above,
which depends on four parameters $0\le A,B,C,D\le n$, 
with $A+B\equiv1\bmod2$ and $B\le C$.
For the number of triangulations of such a hook shape we get a recursion
\begin{eqnarray*}
\lefteqn{h(A,B,C,D) \ \ =\ \  
\binom{\frac{3A+B-1}2}{A} 
\binom{\frac{A+3B-1}2}{B}
\binom{C+D}{C}}\\
&&+\ 
\sum_{0\le a\le A,    \ 0\le b\le B\atop
     a+b\equiv1\bmod2,\ a+b<A+B}
h(a,b,C,D)
\binom{\frac{3A+B-3a-b}{2}-1}{A-a}
\binom{\frac{A+3B-a-3b}{2}-1}{B-b}\\
&&+\ 
\sum_{0\le a\le D,\ 0\le b\le \frac{A+B-1}{2}\atop
     a+b\equiv1\bmod2,\ \frac{a+b+1}{2}\le B}
h(a,b,\tfrac{A+B-1}{2},A)
\binom{D+C-\frac{3a+b+1}{2}}{D-a}
\binom{\tfrac{A+3B-a-3b}{2}-1}{B-\tfrac{a+b+1}{2}}
\\
&&+\ 
h(\tfrac{3B-A-1}{2},\tfrac{A+B-1}{2},\tfrac{A+B-1}{2},A)
\binom{C+D-\frac{5B-A-1}{2}}{C-B}\quad
\textrm{if }D\ge\tfrac{3B-A-1}{2}\ge0.
\end{eqnarray*}
The four terms in this recursion correspond to the 
four cases depicted in the figure below, 
where the fourth case --- of a long diagonal of width~$3$ ---
occurs only in the case where the second endpoint of the
diagonal, which may be computed to have $y$-coordinate
$\frac{3B-A-1}2$, comes to lie within the hook.
\medskip

\noindent%
\input{v_strip3recursion3.pstex_t}%

\subsection{Strips of (fixed) width \boldmath$m$. }

We now describe a recursive strategy for the enumeration of unimodular
triangulations of grids of arbitrary size.  The method is applicable
for triangulations of general finite point sets --- but it is
effective only in the special case where the points lie on a small
family of parallel (vertical, say) lines; in our case this is the
situation of small (fixed) $m$ and variable~$n$. The key observation is
that any triangulation may be dismantled by removing triangles from
the upper boundary, while maintaining a lattice triangulation of a
$y$-convex lattice polygon.  Since for fixed~$m$ the number of such
polygons in~$P_{m,n}$ is bounded by a polynomial in~$n$, this 
yields an efficient dynamical programming algorithm in this
case.

Let $\Delta_1,\Delta_2\subset\R^2$ be two triangles whose intersection
$\Delta_1\cap\Delta_2$ is a common (empty, zero-, or one-dimensional)
face of both~$\Delta_1$ and~$\Delta_2$.  
We say that \emph{$\Delta_2$ lies above $\Delta_1$} ($\Delta_1\prec
\Delta_2$) if there are two points $(x,y_1)\in
\Delta_1\setminus\Delta_2$ and $(x,y_2)\in \Delta_2\setminus\Delta_1$
on a vertical line with $y_1<y_2$.
For example, in our figure the shaded triangle lies above the 
other one; the other two pairs of triangles are incomparable.
\[  \includegraphics[height=3cm]{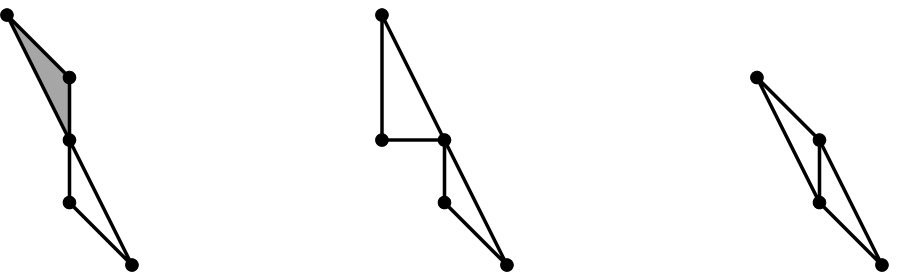}\]
Due to
the convexity of the triangles and the intersection condition imposed
on them, this is a well-defined asymmetric and irreflexive relation.

\begin{lemma}
  \label{lem:abov}
  There is no sequence $\Delta_0,\dots,\Delta_{t-1}\subset\R^2$ of
  triangles (such that the intersection of any two among them is a face
  of both) satisfying
  \begin{equation}
    \label{eq:abov}
    \Delta_0\abov\Delta_1\abov\dots\abov\Delta_{t-1}\abov\Delta_0\enspace.
  \end{equation}
\end{lemma}

\begin{proof}
  Suppose that $\Delta_0,\dots,\Delta_{t-1}\subset\R^2$ is a minimal cyclic
  sequence of triangles (such that the intersection of any two among
  them is a face of both), i.e., it satisfies~(\ref{eq:abov}) (it is
  \emph{cyclic}), but no subsequence of
  $\Delta_0,\dots,\Delta_{t-1}\subset\R^2$ is cyclic.  
  The orthogonal projections $x(\Delta_i)$ to the $x$-axis
  have the following three properties (where all
  indices are taken modulo~$t$): 
  \\(a)~$x(\Delta_i)\cap x(\Delta_{i+1})\not=\emptyset$, 
  \\(b)~$x(\Delta_i)\not\subseteq x(\Delta_{i-1}), x(\Delta_{i+1})$, and 
  \\(c)~$x(\Delta_{i-1})\cap x(\Delta_{i+1})=\emptyset$, 
  \\where (a) follows
  immediately from the definition of~$\abov$ and~(b) as well as~(c)
  are due to the minimality of the cycle.
  
  But (a), (b), and (c) together imply that the intervals
  $x(\Delta_0),\dots,x(\Delta_{t-1})$ ``either run left-to-right or
  right-to-left''; in particular, we have
  $x(\Delta_0)\cap x(\Delta_i)=\emptyset$ for $i\in\{2,\dots,t-1\}$,
  contradicting $\Delta_{t-1}\abov\Delta_0$.
\end{proof}

Of course, both in the definition of~$\abov$ as well as in 
Lemma~\ref{lem:abov} one can replace
``triangle'' by ``compact convex set.'' 

The relation~$\abov$ was defined with respect to parallel projection
here.  If one defines an ordering with respect to central rather than
to parallel projection, then the analog to Lemma~\ref{lem:abov} (for
arbitrary centers of projection) does not hold. For Delaunay
triangulations, however, there is an analog to Lemma~\ref{lem:abov}.
(See De Floriani et al.~\cite{DFNP91} for dimension~$2$, and
Edelsbrunner~\cite{Edelsbrunner92} for arbitrary dimensions.)

Due to Lemma~\ref{lem:abov}, the relation~$\abov$ induces 
a partial order on the set of triangles
in~$\R^2$, which we will also denote by~$\abov$.

A sequence $\mathcal{T}_1,\dots,\mathcal{T}_{2mn}$ of sets of
triangles will be called an \emph{admissible sequence} for $P_{m,n}$
if $\mathcal{T}_1$ is a unimodular triangulation of~$P_{m,n}$, and if,
for each $i=2,3,\dots,2nm$, we have
$\mathcal{T}_i=\mathcal{T}_{i-1}\setminus\{\Delta\}$ for some
$\abov$-maximal triangle~$\Delta$ in~$\mathcal{T}_{i-1}$.  A subset
$S\subset\R^2$ is called an \emph{admissible shape} (of $P_{m,n}$) if
it can be obtained as a union
$S=\bigcup_{\Delta\in\mathcal{T}_i}\Delta$ for an admissible sequence
$\mathcal{T}_1,\dots,\mathcal{T}_{2mn}$ and some
$i\in\{1,\dots,2mn\}$.  Every admissible shape is $y$-convex (i.e.,
its intersection with any vertical line is connected). It is
determined by its \emph{upper boundary segments}, i.e., the sequence
of line segments $[l^{(1)},r^{(1)}],\dots,[l^{(t)},r^{(t)}]$ with
$l^{(1)}_x=0$, $l^{(t)}_x=n$, and $r^{(j-1)}_x=l^{(j)}_x$ for
$j\in\{2,\dots,t\}$, such that, for each point~$p$ in the relative
interior of any of the segments, $p+(0,\varepsilon)\not\in S$ holds
for all $\varepsilon>0$.

Let~$S$ be an admissible shape. We denote by $\mathcal{T}_{\max}(S)$
the set of all $\abov$-maximal unimodular triangles in $S$,
that is, the finite set of all unimodular triangles that could
be $\abov$-maximal in \emph{some} unimodular triangulation of~$S$.
For example, the figure below indicates the 12
$\abov$-maximal triangles of the shaded admissible shape.
\[
\includegraphics[height=5cm]{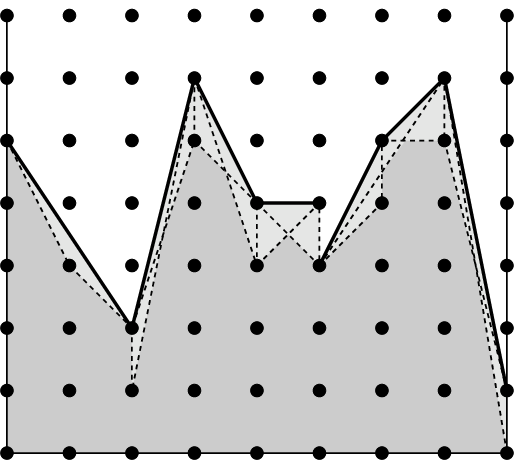}
\]
Any admissible shape~$S'$ that arises
from~$S$ by removing some triangles contained in~$\mathcal{T}_{\max}(S)$ is
called an \emph{admissible subshape} of~$S$. These triangles have
disjoint interiors, and they are uniquely determined for each
admissible subshape~$S'$ of~$S$ (compare the proof of
Lemma~\ref{lem:maxtria}). Their number is denoted by $\#(S',S)$.

Since every unimodular triangulation of~$S$ contains at least one
triangle from~$\mathcal{T}_{\max}(S)$, we obtain the following
inclusion-exclusion formula for the numbers $f(S)$ of unimodular
triangulations of admissible shapes~$S$.

\begin{lemma}
  \label{lem:inex}
  Every admissible shape~$S$ has
  \begin{equation}
    \label{eq:admsh}
    f(S)=\sum_{S'}(-1)^{\#(S',S)-1}f(S')    
  \end{equation}
  unimodular triangulations, where the sum is taken over all
  admissible proper subshapes~$S'$ of~$S$.
\end{lemma}

Lemma~\ref{lem:inex} allows us to compute $f(m,n)=f([0,m]\times[0,n])$
via a dynamic programming approach: Determine the numbers $f(S)$
via~(\ref{eq:admsh}) in some order such that every admissible shape
appears after all its admissible subshapes. In order to analyze the
running time of such an algorithm, we first need to estimate the number
of admissible shapes.

\begin{lemma}
\label{lem:upbd}
  Let $[l^{(1)},r^{(1)}],\dots,[l^{(t)},r^{(t)}]$ be the sequence of upper boundary
  segments of an admissible shape~$S$. Then
  $$
  l^{(j)}_y\in\{r^{(j-1)}_y-1,r^{(j-1)}_y,r^{(j-1)}_y+1\}
  $$
  holds for each $2\leq j\leq t$.
\end{lemma}

\begin{proof}
  This follows by induction on the number of triangles removed in
  order to obtain~$S$: Every vertical edge of a triangle in a
  unimodular triangulation has length one, and thus removing a $\abov$-maximal
  triangle from an admissible shape never creates a vertical
  boundary part of height more than~$1$.
\end{proof}

Lemma~\ref{lem:upbd} implies the following bound on the number of
admissible shapes. 

\begin{lemma}
  \label{lem:numadm}
  There are at most $(3n+2)^{m-1}(n+1)^2$ admissible shapes of~$P_{m,n}$.
\end{lemma}

\begin{proof}
  The upper boundary of an admissible shape has $n+1$ possible start
  and $n+1$ possible end points. At every interior $x$-coordinate
  ($x\in\{1,\dots,m-1\}$) either a segment of the upper boundary ends
  and a new one starts ($3(n+1)-2$ possibilities, by
  Lemma~\ref{lem:upbd}) or a segment ``passes through'' (one possibility).
\end{proof}

The second important quantity for the analysis of the running time of
the dynamic programming algorithm proposed above is the maximal number
of summands that may occur in~(\ref{eq:admsh}).

\begin{lemma}
  \label{lem:maxtria}
  Every admissible shape~$S$ of $P_{m,n}$ has at most 
  $$
  \Big(\frac{3+\sqrt{13}}{2}\Big)^m<\ 3.31^m
  $$
  admissible subshapes.
\end{lemma}

\begin{proof}
  Let $(B_1,\dots,B_t)$ be the sequence (from left to right) of upper
  boundary segments of~$S$.  Each triangle in $\mathcal{T}_{\max}(S)$
  contains at least one of the edges $\{B_1,\dots,B_t\}$.
  
  Let $B\in\{B_1,\dots,B_t\}$ be a segment of the upper boundary
  of~$S$. There are at most two triangles in $\mathcal{T}_{\max}(S)$
  containing~$B$ and one of its adjacent segments.  Each other
  triangle in $\mathcal{T}_{\max}(S)$ that contains $B$ must have its
  third vertex~$v$ in $B+\R\cdot(0,-1)$ on the line that is parallel
  to~$B$ at distance $\ell_2(B)^{-1}$, where $\ell_2(B)$ denotes the
  Euclidean length of~$B$ (because the triangle has area $\frac12$).
  Since~$v$ must be integral and there is no integral point in the
  relative interior of~$B$, there are at most two possibilities
  for~$v$. Let us call one of them \emph{the first triangle
    below~$B$}, and the other one, \emph{the second triangle
    below~$B$} (if they exist).
  
  For every admissible subshape~$T$ of~$S$, define a word
  $w(T)\in\{0,\alpha,\beta,\gamma\}^{\star}$ by replacing $B_i$ by
  \begin{itemize}\itemsep=-2pt
  \item[`$\alpha$' ] if $S\setminus T$ contains the first triangle below~$B_i$,
  \item[`$\beta$' ] if $S\setminus T$ contains the second triangle below~$B_i$,
  \item[`$\gamma$' ] if $S\setminus T$ contains the triangle formed
    by~$B_i$ and~$B_{i-1}$, and
  \item[`$0$'      ] otherwise.
  \end{itemize}
  Clearly, every `$\gamma$' in $w(T)$ has a `$0$' as its left neighbor;
  furthermore, $w(\cdot)$ is an injective mapping. Therefore, the number
  of admissible subshapes of~$S$ is bounded from above by the
  function $\varphi(m)$ defined recursively via
  $$
  \varphi(0)=1\ ,\qquad
  \varphi(1)=3\ ,\qquad
  \varphi(m)=3\varphi(m-1)+\varphi(m-2)
  \quad(m\geq 2).
  $$
  Using standard techniques 
  (see, e.\,g., \cite[Thm. 4.1.1]{StanleyEnumComb}), one derives 
  from this recursion that
  $\varphi(m)<\big(\frac{3+\sqrt{13}}{2}\big)^m  $ for $m\ge1$. 
\end{proof}

Lemmas~\ref{lem:numadm} and~\ref{lem:maxtria} imply that (for
every~$m$) the dynamic programming algorithm needs at most
$$
3.31^m(3n+2)^{m-1}(n+1)^2\ <\ 10^m (n+1)^{m+1}
$$
arithmetic operations. The actual running time of an implementation
heavily depends on the data structures for storing the admissible
shapes and on the way by which one determines the admissible subshapes in
that data structure. Therefore, here we include only the following
rough statement.

\begin{theorem}
  For every fixed~$m$, the function $f(m,n)$ can be computed in time
  bounded by a polynomial in~$n$.
\end{theorem}

In our implementation, we organize the admissible shapes as the leaves
of a tree whose nodes are the prefixes of the sequences of upper
boundaries segments of admissible shapes. The data structure allows
quite efficient access to the admissible subshapes while not wasting
too much memory.  Nevertheless, the bottleneck in the computations
is always memory. It is crucial to use an ordering of the
admissible shapes in which for each shape~$S$ the shapes of which it
is an admissible subshape come as soon as possible after it; this
allows one to keep in memory only a subset of the admissible shapes at
each point of time.

Furthermore, in~(\ref{eq:admsh}) there is no need to sum over \emph{all}
admissible subshapes~$S'$. It suffices to consider those~$S'$ that
arise from removing triangles in~$\mathcal{T}_{\max}(S)$ that are
contained in $\{(x,y)\in\R^2:x\geq c\}$, where~$c$
is the maximal $x$-coordinate where $l^{(j)}_y=r^{(j-1)}_y+1$ for
some~$j$ and the upper boundary segments
$[l^{(1)},r^{(1)}],\dots,[l^{(t)},r^{(t)}]$ of~$S$ (if that maximum exists). 

Of course, when computing the number of unimodular triangulations of
$P_{m,n}$ by our method, we obtain as a byproduct the number of
unimodular triangulations of several interesting polygons inside
$P_{m,n}$, including $P_{m,1}$, \dots, $P_{m,n-1}$. 

The algorithm described above cannot only be used to calculate the
number of unimodular triangulations of any admissible shape~$S$
in~$P_{m,n}$; it can also be extended to produce a uniformly
distributed random triangulation within the same
asymptotic running time, and thus, in polynomial time (depending
on~$n$) for fixed~$m$. For this, one just determines the numbers of
triangulations of those admissible subshapes of~$S$ that arise from
removing one single triangle; with respect to the corresponding random
distribution one then chooses one of the $\abov$-maximal triangles
in~$S$ at random, and proceeds with the subshape obtained by removing
it.

\subsection{Explicit values}
\label{subsec:explval}

We have implemented the algorithms described in
Subsection~\ref{subsec:narrow} in C++, using the \texttt{gmp} library~\cite{gmp}
for exact arithmetic, with the interface to it provided by the
\texttt{polymake} system~\cite{polymake}. 

Results obtained by our code are compiled in the Appendix (see
Tables~\ref{tab:m2}, \ref{tab:m3}, \ref{tab:m4}, \ref{tab:m5},
and~\ref{tab:m6}).  The number of unimodular triangulations of
$P_{m,n}$ asymptotically grows exponentially with $mn$ (see also
Section~\ref{sec:bounds}).  Therefore, it is more convenient to view
the function $f(m,n)$ on a logarithmic scale, normalized by~$mn$.

\defn
The \emph{capacity} of the $m\times n$ grid is 
\[
c(m,n)\ \ :=\ \ \frac{\log_2 f(m,n)}{mn}.
\]
\enddefn

The following Figure shows the capacity functions $c(m,n)$ for
$m\in\{1,\dots,6\}$.  The largest capacity we found is $2.055792$ (for
$m=4$ and $n=32$, see Table~\ref{tab:m4}).

\begin{figure*}[h]
  \centering
  \includegraphics[height=6.8cm]{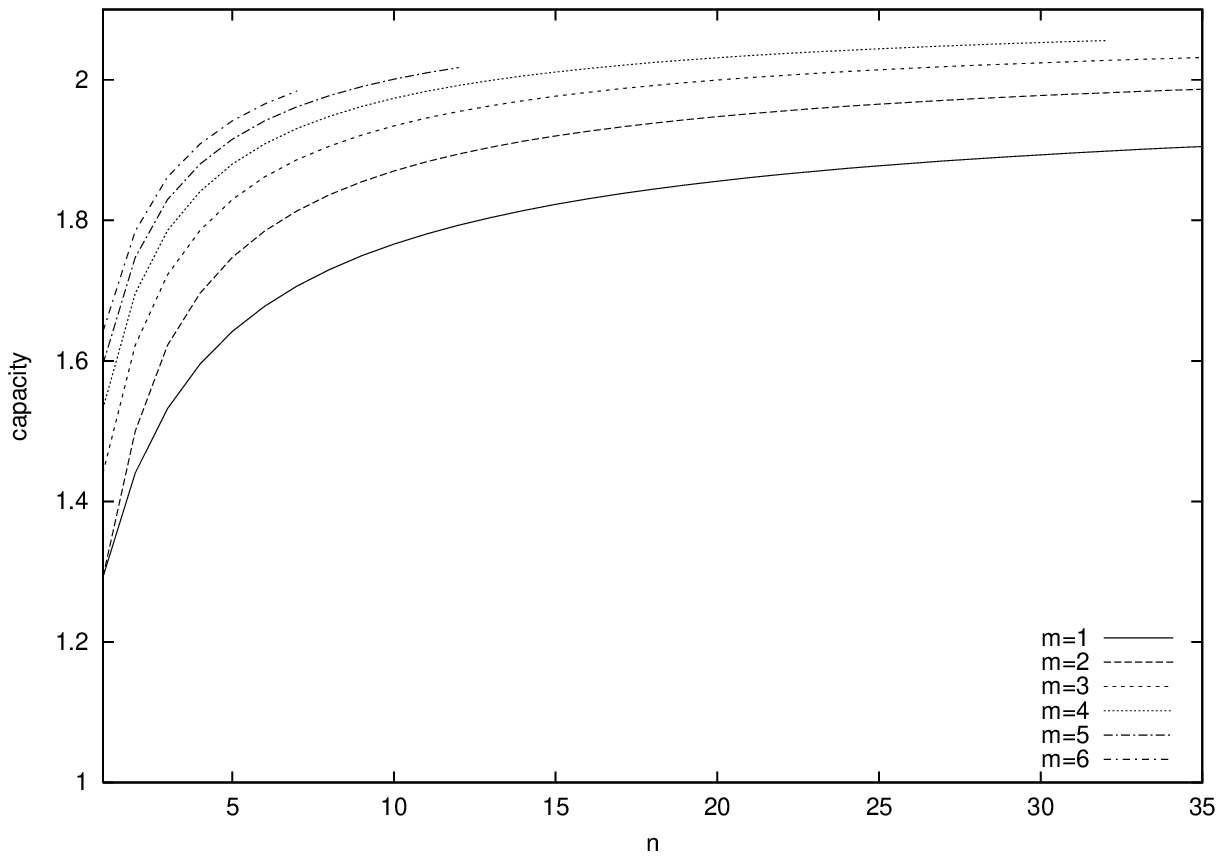}
  \label{fig:cap}
\end{figure*}

To give an impression of the amount of (machine) resources required
for the calculations: The run of the admissible shape algorithm for
$m=6$ and $n=7$ needed about three gigabytes of memory. The grid
$P_{6,7}$ has $370252552$ admissible shapes. We generated them in
lexicographical order with respect to the pairs of starting heights
and volumes. This way, never more than $15\%$ of the admissible shapes
had to be kept in memory simultaneously.  Notice that more than 400
megabytes (of the 3 gigabytes in total) where needed just to store the
\emph{numbers} of triangulations of the admissible shapes in the
memory. The CPU time used for the computation was about 25 hours.  Our
computations were performed on a SUN UltraAX MP machine equipped with
four 448 MHz UltraSPARC-II processors (of which we used only one) and
4 gigabytes main memory.

For very small parameters, Meyer~\cite{Meyer} has enumerated all
unimodular triangulations by Avis and Fukuda's reverse search method
(sketched at the beginning of this section) and checked them for
regularity. Table~\ref{tab:regsmall} shows the results. While for
these small parameters irregular triangulations are quite rare, the
picture changes drastically when~$m$ and~$n$ get larger, see
Section~\ref{sec:trias}. 

\begin{table}[ht]
  \centering
  \begin{tabular}{crrr}
    \hline
    $m\times n$ &
    \multicolumn{1}{c}{$\#$ triangulations} &
    \multicolumn{1}{c}{$\#$ irregular} &
    \multicolumn{1}{c}{fraction} \\
    \hline
    $3\times 3$ & 46456     & 4       & .000086 \\
    $3\times 4$ & 2822648   & 502     & .000178 \\
    $3\times 5$ & 182881520 & 63528   & .000347 \\
    $4\times 4$ & 736983568 & 1553020 & .002107 \\
    \hline
  \end{tabular}
  \caption{Number of regular triangulations of small grids.}
  \label{tab:regsmall}
\end{table}


\section{Bounds}
\label{sec:bounds}

\subsection{Patching}

Any two unimodular triangulations of $P_{m,n_1}$ and $P_{m,n_2}$ can
be patched to a unimodular triangulation of $P_{m,n_1+n_2}$. Thus we
have the follwing supermultiplicativity relation, where $\firreg(m,n)$
denotes the number of irregular unimodular triangulations of $P_{m,n}$. 

\begin{lemma}
  \label{lem:supermul}
  For $m,n_1,n_2\geq 1$ the following relations hold.
  \begin{enumerate}
    \item[\rm(i)] \mbox{\centerline{$f(m,n_1+n_2)\ge f(m,n_1)f(m,n_2)$}}
    \item[\rm(ii)] \mbox{\centerline{$\firreg(m,n_1+n_2)\ge f(m,n_1)\firreg(m,n_2)$}}
  \end{enumerate}
\end{lemma}

With respect to regular triangulations, patching is dangerous, as
demonstrated by the following example of a non-regular triangulation
of $P_{4,4}$ composed of four regular triangulations of $P_{2,2}$
(suggested by Francisco Santos).
\[
\includegraphics[height=33mm]{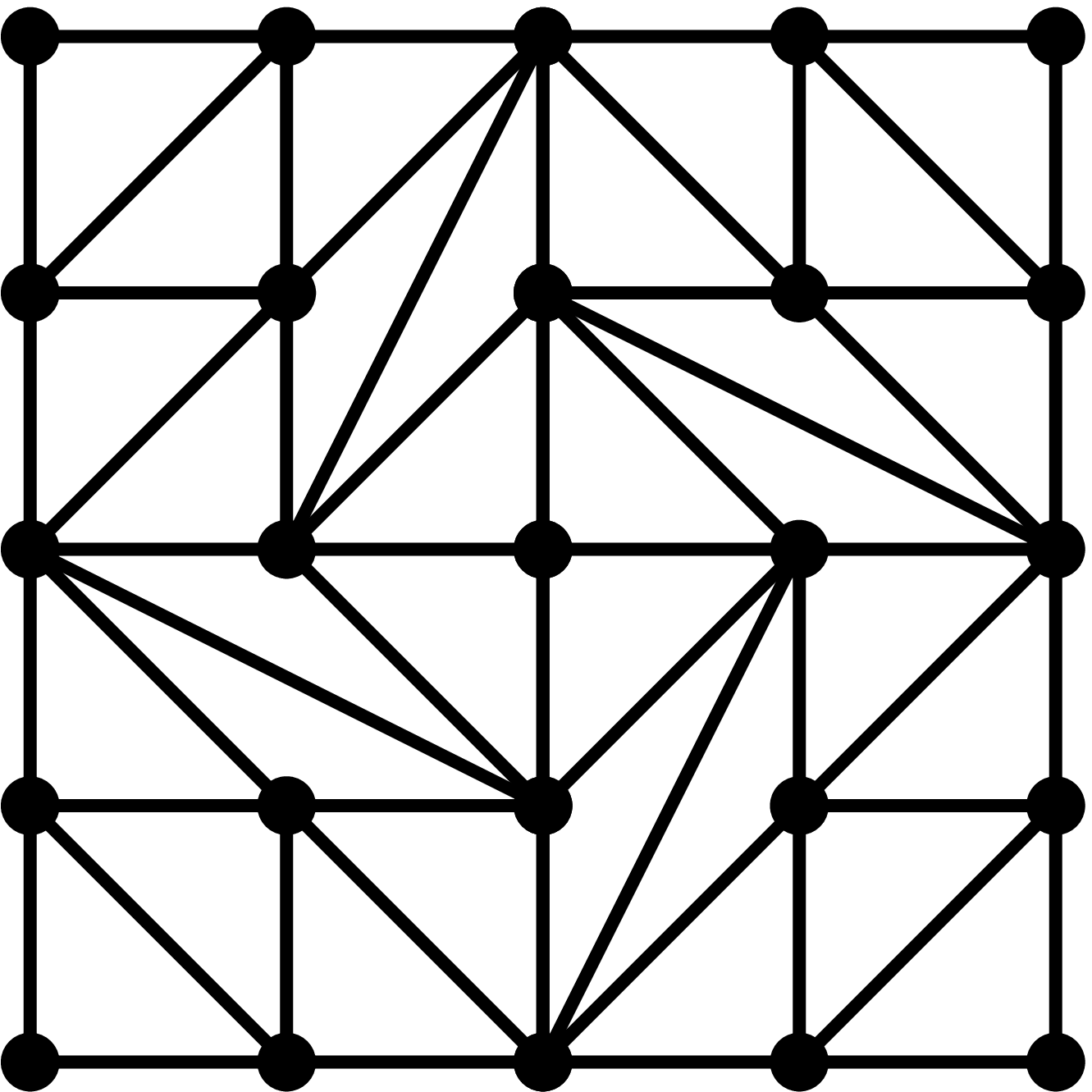}
\]
However, a much more general theorem by Goodman and Pach
\cite{GoodmanPach} says that any two regular
triangulations of two disjoint convex polytopes $P_1,P_2\subset\R^d$
can be extended to a regular triangulation of $\conv(P_1\cup P_2)$
without additional vertices. Thus also for the regular case 
we get a (slightly weaker) supermultiplicativity relation.

\begin{lemma}
  \label{lem:supermulreg}
  For $m,n_1,n_2\geq 1$  we have
  $$
  \freg(m,n_1+n_2+1)\ge \freg(m,n_1)\freg(m,n_2)\ .
  $$
\end{lemma}

The following figure illustrates the patching of
Lemma~\ref{lem:supermulreg}:
\begin{center}
  \includegraphics[height=33mm]{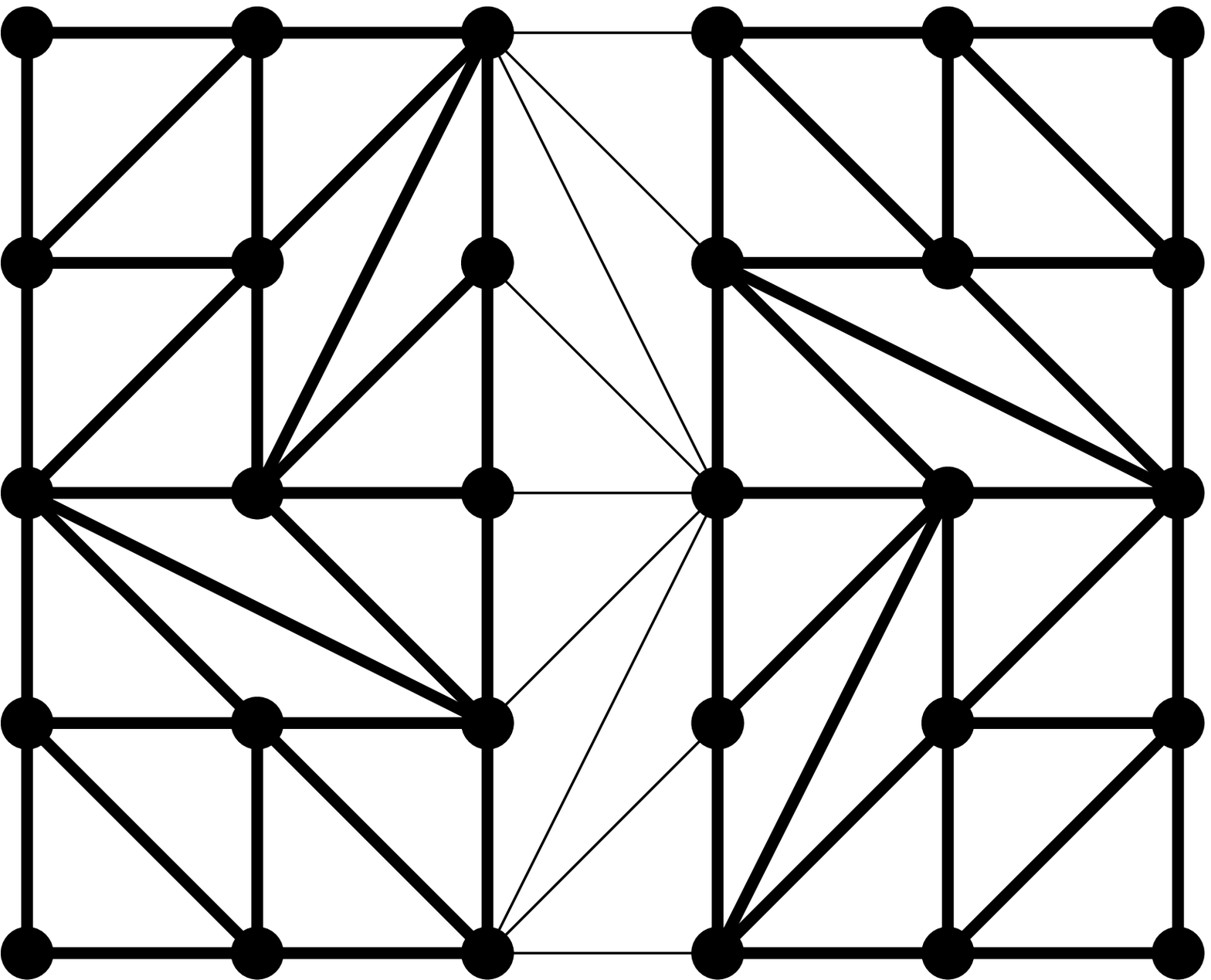}
\end{center}

For $n_2=1$ we will (Lemma~\ref{lem:patchreg1strip}) strengthen the
inequality in Lemma~\ref{lem:supermulreg}.  

Let us fix some notations first.  For any function
$h:P_{m,n}\longrightarrow\R$ we denote by
$H:P_{m,n}\longrightarrow\R^3$ the function with
$H(x,y)=(x,y,h(x,y))$. The function~$h$ is called a \emph{lifting
  function} of a triangulation~$\mathcal{T}$ of~$P_{m,n}$, if
$\mathcal{T}$ is the image of the set of ``lower facets'' of the
$3$-polytope $\conv\{H(x,y):(x,y)\in P_{m,n}\}$ under orthogonal
projection to the $x,y$-plane (deletion of the third coordinate).

A function~$h$ is a lifting function of~$\mathcal{T}$ if and only
if~$h$ is convex and piecewise linear, and its (maximal) domains of
linearity are the triangles in~$\mathcal{T}$. In particular, one may
add to~$h$ any convex piecewise linear function whose domains of linearity
are unions of triangles of~$\mathcal{T}$ in order to obtain another
lifting function for~$\mathcal{T}$.

A triangulation is regular if and only if it has a lifting function.

The following result shows that all unimodular triangulations of a
strip of width $1$ are regular; moreover one has a lot of freedom in
choosing the respective lifting functions, which we will exploit
below.

\begin{lemma}
  \label{lem:stripreg}
  Let $\mathcal{T}$ be a unimodular triangulation of a lattice
  trapezoid with two parallel vertical or horizontal sides $S_0$ and
  $S_1$ at distance one. Every piecewise linear function
  $h_0:S_0\longrightarrow\R$ that is strictly convex on $S_0\cap\Z^2$
  can be extended to a lifting function for~$\mathcal{T}$.
\end{lemma}

\begin{proof}
  Let $p_0,\dots,p_r$ be the integral points on $S_0$, and let
  $\{e_{i,1},\dots,e_{i,k(i)}\}$ be the edges of $\mathcal{T}$
  connecting $p_i$ to $S_1$. Let $S_0^+(i)$ and $S_0^-(i)$ be the
  closures of the two components of $S_0\setminus\{p_i\}$. Then we may
  decompose the function $h_0$ as
  $$
  h_0(x)=\sum_{i=0}^r\sum_{j=1}^{k(i)}h_{i,j}(x)\ ,
  $$
  where each $h_{i,j}$ is a convex function on $S_0$, linear both
  on $S_0^+(i)$ and on~$S_0^-(i)$, and having its  unique break-point
  at~$p_i$.

  Now we extend each $h_{i,j}$ to a convex, piecewise-linear function
  defined on the entire trapezoid such that
  it has a break-line at the edge $e_{i,j}$ and is linear above
  and below this line.
  Then the sum of all $h_{i,j}$ is a lifting function for~$\mathcal{T}$.
\end{proof}

\begin{prop}
\label{prop:fregm12}
  For $n\geq 1$ the following relations hold.
  \begin{enumerate}
  \item[\rm(i)] \mbox{\centerline{$\freg(1,n)=f(1,n)=\binom{2n}{n}$}}
  \item[\rm(ii)] \mbox{\centerline{$\freg(2,n)=f(2,n)$}}
  \end{enumerate}
\end{prop}

\begin{proof}
  Part~(i) follows immediately from Lemma~\ref{lem:stripreg}(i) (and
  equation~(\ref{eq:1n})). 
  
  Since patching two regular triangulations along a \emph{single} edge
  preserves regularity, it suffices for the proof of part~(ii) to show
  that every unimodular triangulation of shapes of one of the  forms
  \begin{center}
    \includegraphics[height=40mm]{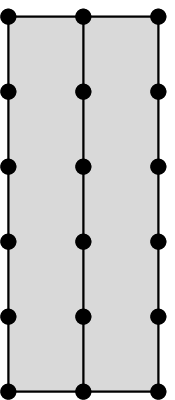}\hspace{20mm}
    \includegraphics[height=40mm]{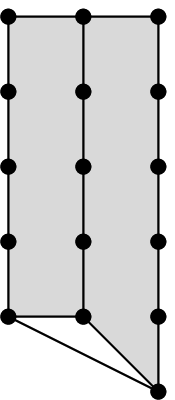}\hspace{20mm}
    \includegraphics[height=40mm]{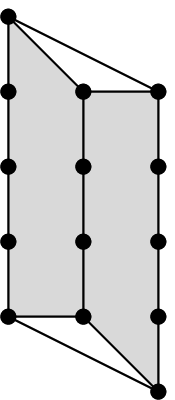}
  \end{center}
  is regular. But this can be derived from
  Lemma~\ref{lem:stripreg}.  One starts from an arbitrary
  prescribed strictly convex function on the middle column, and after
  the extensions to the two shaded vertical trapezoids of width one
  obtained from Lemma~\ref{lem:stripreg}, one adds a piecewise linear
  function that is constant on the left strip, linear on the right
  strip, and sufficiently large on the right column of vertices.
\end{proof}

Similarly, one proves the following  strengthening of
Lemma~\ref{lem:supermulreg} for $n_2=1$, as announced earlier.

\begin{lemma}
  \label{lem:patchreg1strip}
  For $m,n\geq 1$, we have
  $$
  \freg(m,n+1)\geq\freg(m,n)\cdot\binom{2n}{n}\ .
  $$
\end{lemma}

\subsection{Limit capacities}

In the following, we will show that the capacities $c(m,n)$ and
$\creg(m,n)$ (with $f(m,n)=2^{c(m,n)mn}$ and
$\freg(m,n)=2^{\creg(m,n)mn}$) asymptotically behave well, which
allows us to focus on their limits subsequently. Note that all
capacities are bounded (see Theorem~\ref{thm:anclin}).

\begin{prop}
  Let $m\geq 1$.
  \begin{enumerate}
    \item[\rm(i)] The limit \ 
       $c_m:=\displaystyle\lim_{n\rightarrow\infty}  c   (m,n)$ \ exists.
    \item[\rm(ii)] The limit \ 
   $\creg_m:=\displaystyle\lim_{n\rightarrow\infty} \creg(m,n)$ \ exists.
  \end{enumerate}
\end{prop}

\begin{proof}
  Lemmas~\ref{lem:supermul}(i) and~\ref{lem:supermulreg} imply by
  Fekete's lemma \cite[Lemma 11.6]{LW} that
  $$
  \lim_{n\rightarrow\infty}f(m,n)^{\frac{1}{n}}
  \qquad\text{and}\qquad
  \lim_{n\rightarrow\infty}\freg(m,n-1)^{\frac{1}{n}}  
  $$
  exist. Therefore,  
  $$
    \lim_{n\rightarrow\infty}\frac{1}{m}\log f(m,n)^{\frac{1}{n}}
   =\lim_{n\rightarrow\infty}\frac{\log f(m,n)}{mn}
   =\lim_{n\rightarrow\infty}c(m,n)
  $$
  and
  $$
    \lim_{n\rightarrow\infty}\frac{n}{m(n-1)}\log \freg(m,n-1)^{\frac{1}{n}}
   =\lim_{n\rightarrow\infty}\frac{\log \freg(m,n-1)}{m(n-1)}
   =\lim_{n\rightarrow\infty}\creg(m,n)
   $$
  exist as well.
\end{proof}

While the last proposition concerned the asymptotics of growing~$n$ for
fixed~$m$, the next result shows that also growing~$m$ and~$n$
simultaneously yields nice asymptotics.

\begin{prop}
  Let $m\geq 1$.
  \begin{enumerate}
  \item[\rm(i)] The limit
    $c:=\displaystyle\lim_{m\rightarrow\infty}c(m,m)$ exists. It
    satisfies
    $$
    c=\displaystyle\lim_{m\rightarrow\infty}c_m
    \qquad\text{and}\qquad
    c_{m_0}\leq c\quad (m_0\in\N)\ .
    $$
  \item[\rm(ii)] The limit
    $\creg:=\displaystyle\lim_{m\rightarrow\infty}\creg(m,m)$
    exists. It satisfies
    $$
    \creg=\displaystyle\lim_{m\rightarrow\infty}\creg_m
    \qquad\text{and}\qquad
    \creg_{m_0}\leq \creg\quad (m_0\in\N)\ .
    $$

  \end{enumerate}
\end{prop}

\begin{proof}
  From Lemma~\ref{lem:supermul}(i) one derives (for $m_0,n_0\geq 1$)
  the inequality
  \begin{equation}
    \label{eq:diaglim1}
    f(m,n)\geq f(m_0,n_0)^{\lfloor\frac{m}{m_0}\rfloor\lfloor\frac{n}{n_0}\rfloor}.
  \end{equation}
  For integers $p,q\geq 1$ we define $\Phi(p,q):=1-\frac{p\mod q}{p}$. We have
  $\Phi(q,q)=1$ and $\lim_{p\rightarrow\infty}\Phi(p,q)=1$ for all
  $q\in\N$. 

  Equation(\ref{eq:diaglim1}) then implies
  \begin{equation}
    \label{eq:diaglim2}
    c(m,n)\geq\Phi(m,m_0)\Phi(n,n_0)c(m_0,n_0)\ ,
  \end{equation}
  and, in particular,
  \begin{equation}
    \label{eq:diaglim3}
    c(m_0,n)\geq\Phi(n,n_0)c(m_0,n_0)\ .
  \end{equation}
  Inequality~(\ref{eq:diaglim3}) (together with
  $\displaystyle\lim_{n\rightarrow\infty}\Phi(n,n_0)=1$) yields
  \begin{equation}
    \label{eq:diaglim4}
    c_m\geq c(m,n_0)\ .
  \end{equation}
  Inequality~(\ref{eq:diaglim2}) implies (together
  with $\displaystyle\lim_{m\rightarrow\infty}\Phi(m,m_0)\Phi(m,n_0)=1$)
  $$
  \liminf_{m\rightarrow\infty} c(m,m)\geq c(m_0,n_0)\ ,
  $$
  and therefore,
  \begin{equation}
    \label{eq:diaglim5}
    \liminf_{m\rightarrow\infty} c(m,m)\geq c_{m_0}\ .
  \end{equation}
  
  Finally, we obtain the following chain of inequalities, which, 
  together with (\ref{eq:diaglim5}) proves part~(i) of the proposition.
  The middle inequality is from~(\ref{eq:diaglim5}), and the
  outer ones are due to~(\ref{eq:diaglim4}).
  $$
          \liminf_{m\rightarrow\infty}c_m
    \geq  \liminf_{m\rightarrow\infty}c(m,m)
    \geq  \limsup_{m\rightarrow\infty}c_m
    \geq  \limsup_{m\rightarrow\infty}c(m,m)
  $$

  Part~(ii) is proved similarly, starting from Lemma~\ref{lem:supermulreg}.
\end{proof}

Note that a similar proof yields
$$
c=\lim_{m\rightarrow\infty}c(\alpha m,\beta n)
\qquad\text{and}\qquad
\creg=\lim_{m\rightarrow\infty}\creg(\alpha m,\beta n)
$$
for each pair $\alpha,\beta>0$.

By Lemma~\ref{lem:supermul}(ii), the corresponding ``irregular limit
capacities'' do exist as well, and they are equal to $c_m$
($m\ge3$) and $c$, respectively. Therefore, we do not treat them
explicitly.

\subsection{Lower bounds}\label{subsec:lower_bounds}

\begin{prop}
  The following estimates hold.
  \begin{enumerate}
    \item[\rm(i)] $c_1=\creg_1=2$
    \item[\rm(ii)] $c_2=\creg_2>2.044$
    \item[\rm(iii)] $c_3>2.051$, $c_4>2.055$
    \item[\rm(iv)] $c_m>2.048$ for $m\geq 5$
    \item[\rm(v)] $c_m\ge\creg_m>2$ for $m\geq 3$
  \end{enumerate}
\end{prop}

\begin{proof}
  Part~(i) follows from 
  $$
  f(1,n)=\freg(1,n)=\binom{2n}{n}\approx\frac{2^{2n}}{\sqrt{2n}}\ .
  $$
  Parts~(ii) and~(iii) are results of the computer calculations
  reported in Section~\ref{sec:values}, combined with
  Proposition~\ref{prop:fregm12}(ii).

  Lemma~\ref{lem:supermul}(i) applied to the ``transposed grids'' implies 
  \begin{equation}
    \label{eq:lowerbounds}
    c(m_1+m_2,n)\geq\frac{m_1}{m_1+m_2}c(m_1,n)+\frac{m_2}{m_1+m_2}c(m_2,n)\ ,    
  \end{equation}
  and thus 
  $$
  c_{m_1+m_2}\geq\frac{m_1}{m_1+m_2}c_{m_1}+\frac{m_2}{m_1+m_2}c_{m_2}\ .
  $$
  With~(ii) and~(iii), this yields~(iv).

  Similarly, Lemma~\ref{lem:patchreg1strip} leads to
  $$
  \creg_{m+1}\geq \frac{m}{m+1}\creg_m+\frac{1}{m+1}\cdot 2\ .
  $$
  With~(ii) this proves part~(v).
\end{proof}

Equation~(\ref{eq:lowerbounds}) implies that $c_{kn}\ge c_n$ for
$k\ge2$; for example, this implies that $c_4\ge c_2$ --- but it is not
obvious that $c_3\ge c_2$. Even stronger, one would assume that
\[
2 = c_1 <  c_2 <  c_3 < \ \cdots\ \ \le c\ ,
\]
and
\[
2 = \creg_1 <  \creg_2 <  \creg_3 < \ \cdots\ \ \le \creg\ ,
\]
but neither monotonicity is proved.

\subsection{Upper bounds}\label{subsec:upper_bounds}

From general principles (see Ajtai et al. \cite{ACNS}) one
gets that the capacity for any configuration of $N$ points in the plane
is finite. In the general position case (no three points on a line)
the currently best upper bound is $o(2^{59N})$,
due to Santos \& Seidel \cite{SantosSeidel}.
However, in the very ``degenerate'' case of lattice triangulations,
there are far fewer triangulations:
Orevkov \cite{Orevkov} obtained the bound $f(m,n)\le4^{3mn}=2^{6mn}$.
Very recently, this has been substantially improved by
Anclin \cite{Anclin}, as follows.

\begin{theorem}[Anclin \cite{Anclin}]
\label{thm:anclin}
For all $m,n\ge1$,
\[
f(m,n)\ <\ 2^{3mn-m-n}.
\] 
\end{theorem}

\begin{proof}
Our sketch relies on the essential ideas of Anclin's proof.

The first, crucial observation is that the
midpoints of the edges in any unimodular triangulation 
are exactly the half-integral, not integral points
in~$\conv(P_{m,n})$. (Clearly the midpoint of every edge
is half-integral; the converse
may be derived from Pick's theorem, or from the fact that all 
unimodular triangles are
equivalent to~$\conv\{(0,0),(1,0),(0,1)\}$.)
The number of these half-integral points in the interior of~$P_{m,n}$ 
is $e=3mn-m-n$.

Now any triangulation is built as follows:
The half-integral points are processed in 
an order that is given by a parallel \emph{sweep}.
(See de Berg et al.~\cite[Sect.~2.1]{4Marks}
for a discussion and many further sweeping applications
of this fundamental technique.)

Whenever a point is processed, we add to the partial
triangulation a new edge with the given midpoint.
The key \emph{claim} is that at each such step,
when a half-integral point $v$ is processed, 
there are (at least one and) at most two possibilities
for the new edge with midpoint $v$ to be added.
If this claim is true, then the number of triangulations
is bounded by $2^e$.
\[
\input{anclin1.pstex_t}\qquad\input{anclin2.pstex_t}
\]
To prove the claim, one can verify the following:
Let $[v',\bar v']$ and $[v'',\bar v'']$
be two potential edges with midpoint $v$ that could be added,
where $v'$ and $v''$ are the endpoints below the sweep-line $\ell$;
let $Q$ be the convex hull of the integral points in
the triangle $[v,v',v'']$; then all the $k\ge2$
vertices $v'=v_0,v_1,\dots,v_{k-1},v_k=v''$
of $Q$ are ``visible'' from $v$, and any one of the edges $[v_i,\bar v_i]$
with midpoint $v$ could potentially be added at this step.
Furthermore, the midpoints $\frac12(v_{i-1}+v_i)$ lie below
the sweep-line, and thus the edges $[v_{i-1}v_i]$
are present in the current partial triangulation.
\[
\input{anclin3.pstex_t}
\]
Now assume that $[v_0,\bar v_0]$, $[v_1,\bar v_1]$, $[v_2,\bar v_2]$
are three \emph{adjacent} edges with midpoint $v$ that could be added when
processing $v$. By central symmetry with respect to~$v$
we may assume that the midpoint $w$
of $[v_1,\bar v_2]$ lies below the sweep-line. But the
triangle $[v_1,v,\bar v_2]$ is then an empty triangle of area
$\frac14$, just like the triangle $[v_1,v,v_2]$. From this
we conclude that in the current partial triangulation, 
the edge $[v_1,\bar v_2]$ must be present --- which creates
a crossing with the potential edge $[v_0,\bar v_0]$, and thus 
a contradiction.
\end{proof}

The following upper bounds on the limit capacities follow immediately
from Theorem~\ref{thm:anclin}.

\begin{cor}
For all $m,n\ge 1$, the following inequalities hold:
\begin{enumerate}
  \item[\rm (i) ] $\creg(m,n)\ \le\ c(m,n)\ \le\ 3-\frac1m -\frac1n$
  \item[\rm (ii) ] $\creg_m\ \le\ c_m\ \le\ 3-\frac1m$ \quad 
      {\rm  (in particular, $\creg_2\le c_2\le 2.5$)}
  \item[\rm (iii) ] $\creg\ \le\ c\ \le\ 3$
\end{enumerate}
\end{cor}

As Anclin noted, his proof works much more generally:
For any partial triangulation of a not necessarily simple or convex
lattice polygon, the number of completions is at most $2^{e'}$,
where $e'$ is 
the number of edges that are to be added.





\section{Explicit Triangulations}
\label{sec:trias}

For small grids, one can enumerate all unimodular triangulations by
the reverse search algorithm 
sketched in Section~\ref{sec:values}. For
larger grids, it is desirable to obtain from the huge set of
unimodular triangulations ``random'' ones. There are ways to produce
them, however, in most cases the probability distribution from which
they are chosen is unknown.

\subsection{Generating Random Triangulations}
\label{subsec:rw}
A standard way to compute complex random objects such as triangulations
is to set up a random walk. In our case of unimodular triangulations
of $P_{m,n}$ the method is described easily: First one determines any
starting triangulation~$\mathcal{T}$ of $P_{m,n}$. Then, the following
operation is performed $\tau$ times: Choose one of the (inner) edges
of~$\mathcal{T}$ uniformly at random; if this edge is flippable (see
Section~\ref{sec:values}), then with probability $\frac12$ the current
triangulation $\mathcal{T}$ is replaced by the one obtained
from it by flipping that edge.

As the flip graph of the triangulations is connected (see
Section~\ref{sec:values}), it follows from general principles that,
with $\tau$ tending to infinity, the probability distribution defined by
the output of this algorithm converges to the uniform distribution
(with respect to the ``total variation distance''); see, e.g., Jerrum
and Sinclair~\cite{JerrumSinclair} or Behrends~\cite{Behrends-MarkovChains}.
Unfortunately, not much is known about the speed of convergence of the
distribution. In particular, for general~$m$ and~$n$ it is not known
whether there is a polynomial bound (in $n+m+\varepsilon^{-1}$) on the number~$\tau$
of steps needed to guarantee that the total variation distance between
the produced and the uniform distribution is at most~$\varepsilon$ (i.e.,
if the associated Markov chain is \emph{rapidly mixing}).  The only
exception is the case~$m=1$: Here it follows from results of
Felsner and Wernisch~\cite{FelsnerWernisch} that the Markov chain is
indeed rapidly mixing.

Despite this lack of knowledge on the distribution of the output of
the  random walk algorithm, one still can use it in order to produce
examples of ``interesting'' triangulations.

\subsection{Empirical Results}

For each $n\in\{10,\dots,20\}$, Meyer~\cite{Meyer} generated~$1000$
random unimodular triangulations of $P_{n,n}$ by running the random
walk described in Subsection~\ref{subsec:rw} $10^9$ steps, recording
every $10^6$th triangulation; see Table~\ref{tab:large}.
\begin{table}[ht]
  \centering
  \begin{tabular}{cccc}
    \hline
    $m\times n$ &
    \multicolumn{1}{c}{irregularity} & 
    \multicolumn{1}{c}{max. edge length} & 
    \multicolumn{1}{c}{av. edge length} \\
    \hline
    $10\times 10$ & .355 & 5.538 & 1.614 \\
    $11\times 11$ & .435 & 5.843 & 1.630 \\
    $12\times 12$ & .559 & 6.118 & 1.645 \\
    $13\times 13$ & .696 & 6.397 & 1.659 \\
    $14\times 14$ & .782 & 6.650 & 1.670 \\
    $15\times 15$ & .875 & 6.911 & 1.681 \\
    $16\times 16$ & .927 & 7.151 & 1.690 \\
    $17\times 17$ & .965 & 7.391 & 1.700 \\
    $18\times 18$ & .971 & 7.618 & 1.708 \\
    $19\times 19$ & .992 & 7.821 & 1.713 \\
    $20\times 20$ & .997 & 8.060 & 1.723 \\
    \hline
  \end{tabular}
  \caption{Results for random unimodular triangulations of large grids. The first
  column shows the (empirical) probability of irregularity, the
  second and third columns contain the (empirical) expected values of
  the maximal and the average edge length.}
  \label{tab:large}
\end{table}

The results  support the conjecture that, for~$n$ tending to
infinity, a random unimodular triangulation of $P_{n,n}$ is irregular
with probablity one. A second observation is that the (expected value)
of the average length of an edge in a random triangulation seems to
grow very slowly with~$n$. 

\subsection{Obstructions to Regularity}

All figures shown below have been produced by Meyer~\cite{Meyer}, who
implemented procedures for checking regularity by solving linear
programs using the \cplex\ library.

Proposition~\ref{prop:fregm12} shows that $P_{2,n}$ has only regular
triangulations.  The grid $P_{3,3}$ has precisely the following four
(pairwise congruent) irregular unimodular triangulations.
\begin{center}
  \includegraphics[width=.9\textwidth]{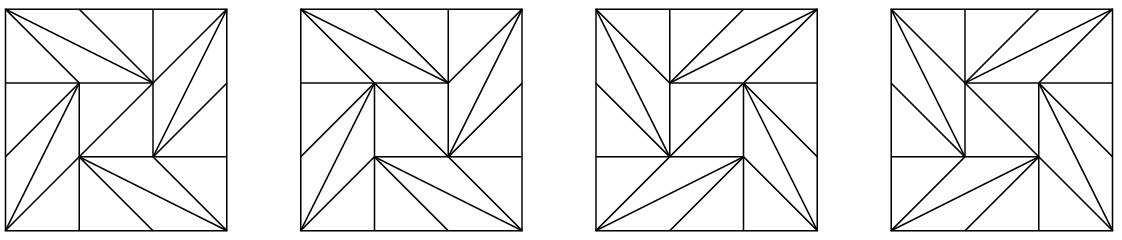}
\end{center}

When trying to understand the reasons for irregularity, it seems
useful to consider (smallest) forbidden patterns for regular
triangulations.  Let~$\mathcal{T}$ be a set of unimodular triangles of
$P_{m,n}$ (such that any two of them intersect in a common face).  We
denote by $S\subset P_{m,n}$ the set of all grid points covered by triangles
in~$\mathcal{T}$. The set~$\mathcal{T}$ is called \emph{regular} if
there is a height function $h:S\longrightarrow\R$ such that for each
triangle $\Delta\in\mathcal{T}$, all $h$-lifted points in 
$S{\setminus}\Delta$ lie strictly above the affine hull of the $h$-lifting
of~$\Delta$.  A subset of~$\mathcal{T}$ is called a 
\emph{minimal irregular configuration} if it is not regular, but all its proper
subsets are regular. Clearly, a regular unimodular triangulation of
$P_{m,n}$ cannot contain any minimal irregular configuration.

Examples for minimal irregular configurations of $P_{3,3}$ are the
following:
\begin{center}
  \includegraphics[width=.9\textwidth]{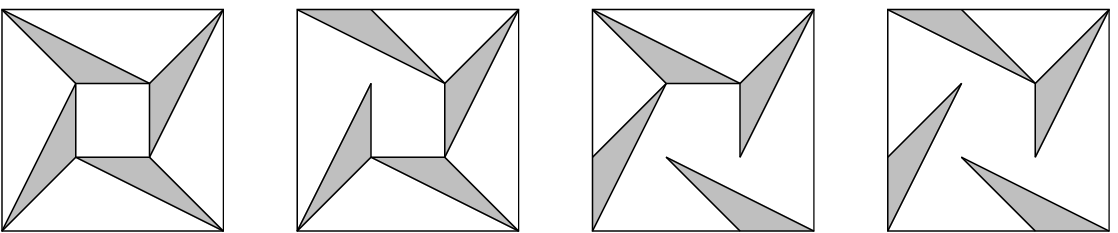}
\end{center}

Already for $P_{3,4}$ many other minimal irregular
configurations occur; some are depicted here:
\begin{center}
  \includegraphics[width=.9\textwidth]{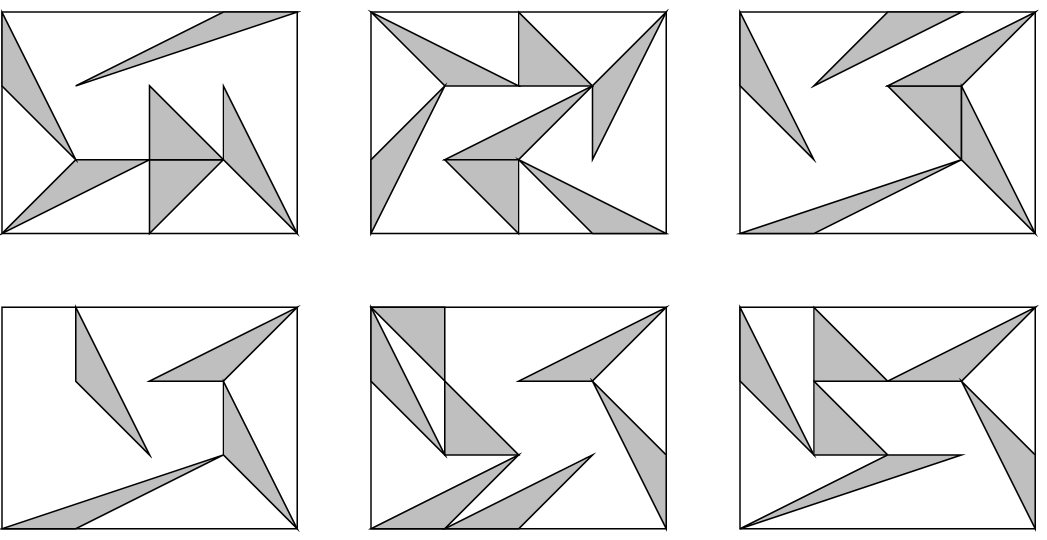}
\end{center}

While these figures still have some similarities with the nice
``whirlpools''
ones for $P_{3,3}$, the picture gets more and more complicated with
growing grid sizes, as the following examples demonstrate:

\begin{center}
  \includegraphics[width=.4\textwidth,bb=0 25 204 145,clip]{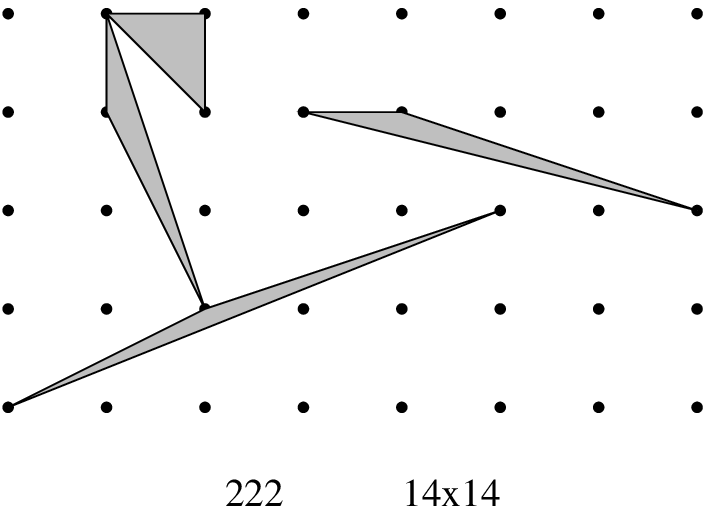}\hspace{.1\textwidth}%
  \includegraphics[width=.4\textwidth,bb=0 25 204 173,clip]{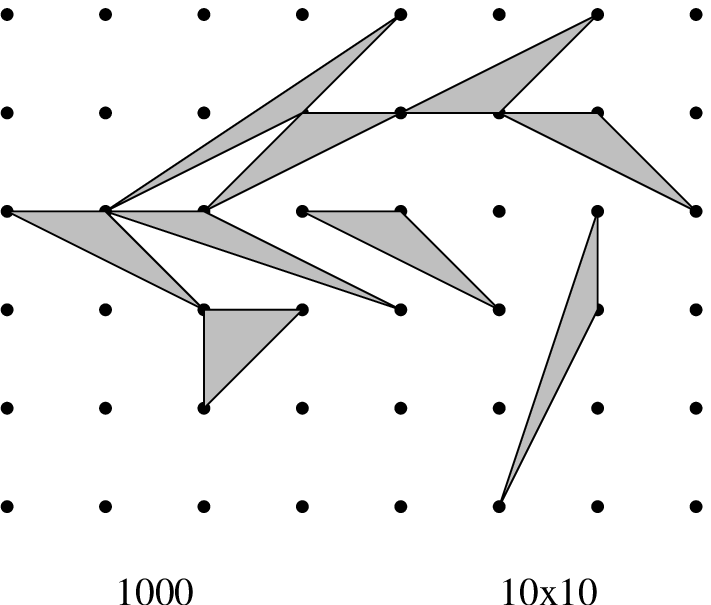}  
\end{center}

\medskip
\begin{center}
  \includegraphics[width=.5\textwidth,bb=0 25 318 259,clip]{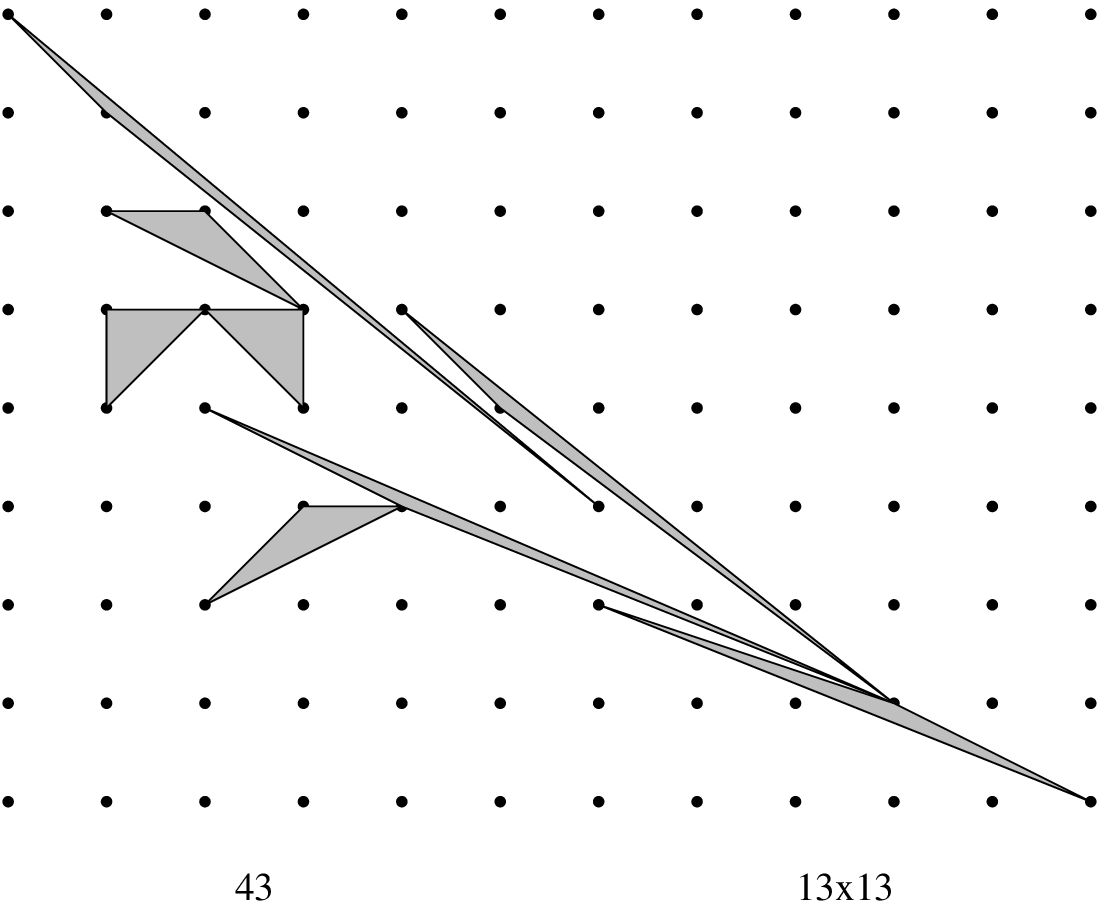}
\end{center}

Viewing these figures, it seems unlikely that one can find any compact
characterization of regularity for unimodular triangulations of
$P_{m,n}$ in terms of forbidden substructures.

We close our zoo of ``explicit triangulations'' with some pairs of
triangulations, found by Meyer's implementation of the random walk.
In each of the figures, the left triangulation is regular, but the right one is
not, although it can be obtained from the left one by flipping just
one edge (drawn bold in the upper left and in the lower right corner,
respectively). For both irregular triangulations, a minimal irregular
configuration contained in it is depicted as well.

\begin{center}
  \includegraphics[width=.9\textwidth]{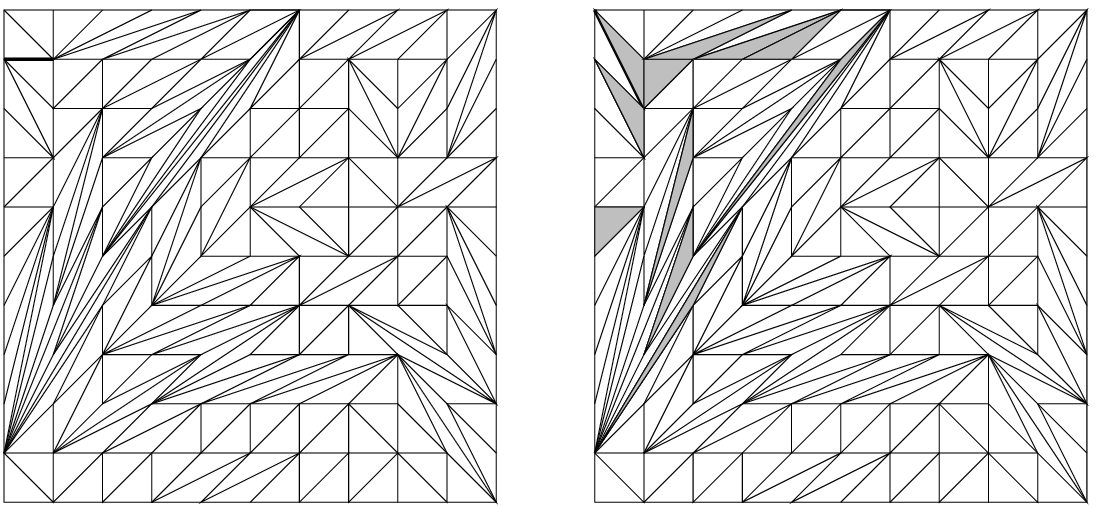}
\end{center}
\begin{center}
  \includegraphics[width=.9\textwidth]{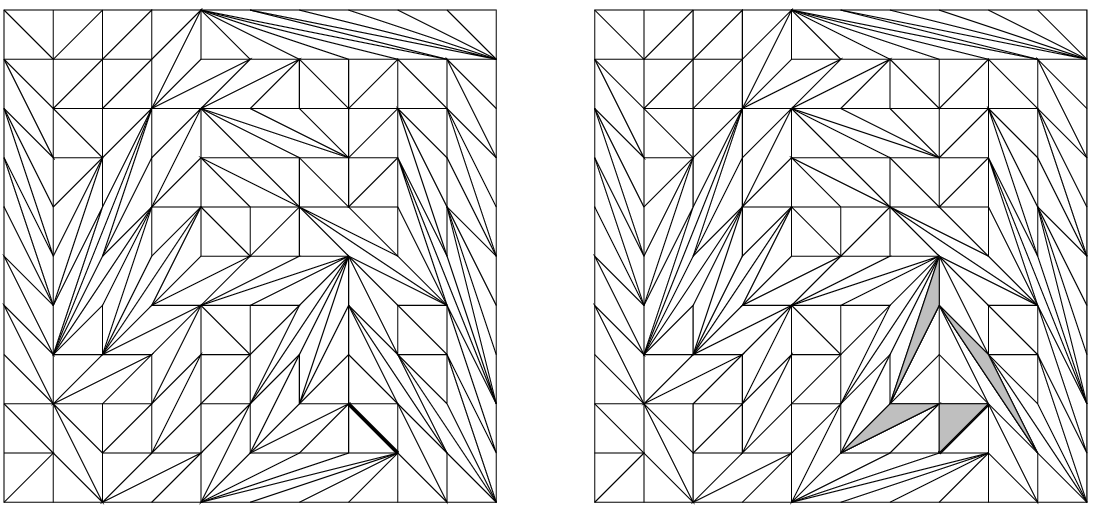}
\end{center}


\subsection*{Acknowledgements.}
We are grateful to 
Oswin Aichholzer,
\'Emile Anclin, 
Stephan Meyer, and
J\"org Rambau
for valuable discussions. Thanks to Francisco Santos and to the
referee for helpful comments
on an earlier version of the manuscript.


\myaddress
\newpage

\section*{Appendix}
\label{sec:appendix}

Below we report on some values $f(m,n)$ and $c(m,n)$ that we have
computed by the algorithms for narrow strips described in
Subsection~\ref{subsec:narrow}. Note that Aichholzer's results referred to in
the captions of the tables have been obtained by a code that works for
\emph{general} point sets in the plane.

\medskip

\begin{table}[h]
  \centering
  {\footnotesize
    \begin{tabular}{rp{.75\textwidth}p{.12\textwidth}@{}}
      \hline
      \multicolumn{1}{c}{$n$\hspace{10mm}} &
      \multicolumn{1}{c}{$\#$ unimodular triangulations  of $P_{2,n}$} &
      \multicolumn{1}{c}{capacity} \\
      \hline
      2 & 64 & 1.500000 \\
3 & 852 & 1.622451 \\
4 & 12170 & 1.696380 \\
5 & 182132 & 1.747462 \\
6 & 2801708 & 1.784822 \\
7 & 43936824 & 1.813494 \\
8 & 698607816 & 1.836244 \\
9 & 11224598424 & 1.854774 \\
10 & 181815529916 & 1.870184 \\
11 & 2964167665340 & 1.883216 \\
12 & 48580814410080 & 1.894393 \\
13 & 799696199314500 & 1.904094 \\
14 & 13212398835196240 & 1.912597 \\
15 & 218976668040908248 & 1.920118 \\
16 & 3639020246503687098 & 1.926820 \\
17 & 60616163842958990268 & 1.932833 \\
18 & 1011775545312594580868 & 1.938260 \\
19 & 16918718677672553292440 & 1.943185 \\
20 & 283368129709983000763876 & 1.947675 \\
21 & 4752924784523774226889308 & 1.951787 \\
22 & 79824154012907603962950312 & 1.955568 \\
23 & 1342199498257069824064033644 & 1.959057 \\
24 & 22592402326314503187343665228 & 1.962288 \\
25 & 380653341141186360494812030908 & 1.965287 \\
\multicolumn{1}{c}{\vdots} & \multicolumn{1}{c}{\vdots} &
\multicolumn{1}{c}{\vdots} \\
\\
375 & 
\parbox{.7\textwidth}{%
322379469795023850474295003353206721961082558051423084880244\\
\ 622493239840267086773982523966267988496124181049802588221265\\
\ 093372321162311900056540572387505992321869061059550930460251\\
\ 208699393384897025341791259853666381631934919090333735984559\\
\ 109034657337254608555310086159512051235323132000393390812675\\
\ 635858586268866021116515894234547255346166243854345190625384\\
\ 072567100646447104480055590853368584565574390095154916379988\\
\ 529322533797949278877762995003200075959780
}\vspace{.1cm}
& 2.044130 \\

      \hline
    \end{tabular}
  }
  \caption{Results for $m=2$ (up to $n=15$ 
     by Aichholzer~\cite{Aichholzer-WWW}).)}
  \label{tab:m2}
\end{table}

\begin{table}[h]
  \centering
  {\footnotesize 
    \begin{tabular}{rp{.75\textwidth}p{.12\textwidth}@{}}
      \hline
      \multicolumn{1}{c}{$n$\hspace{10mm}} &
      \multicolumn{1}{c}{$\#$ unimodular triangulations of $P_{3,n}$} &
      \multicolumn{1}{c}{capacity} \\
      \hline
      3 & 46456 & 1.722619 \\
4 & 2822648 & 1.785718 \\
5 & 182881520 & 1.829755 \\
6 & 12244184472 & 1.861743 \\
7 & 839660660268 & 1.886238 \\
8 & 58591381296256 & 1.905656 \\
9 & 4140106747178292 & 1.921429 \\
10 & 295372308876234428 & 1.934510 \\
11 & 21234538315776214604 & 1.945546 \\
12 & 1535939689343151109944 & 1.954989 \\
13 & 111655493479477379881272 & 1.963164 \\
14 & 8150727077307189203809876 & 1.970314 \\
15 & 597087996550303632801161860 & 1.976623 \\
16 & 43871350204895836758556369212 & 1.982234 \\
17 & 3231797978935266793268797809260 & 1.987258 \\
18 & 238606105193380387756570932194588 & 1.991783 \\
19 & 17651135152017098450035730535703808 & 1.995882 \\
20 & 1308029292984065630362694842042395056 & 1.999613 \\
21 & 97080539975603502667567153853690549804 & 2.003024 \\
22 & 7215158047881650609075575773153609553148 & 2.006154 \\
23 & 536905685776901371485436849505792415847140 & 2.009039 \\
24 & 39997858254082097021224132017959794867440460 & 2.011705 \\
25 & 2982752306685557862989393328648927558138800612 & 2.014178 \\
26 & 222638546950211814977693932477091620801626551100 & 2.016478 \\
27 & 16632293481947394846909242053460530217349594447732 & 2.018623 \\
28 & 1243490745851056260557782821562507156418923407432920 & 2.020627 \\
29 & 93034737749193459157244717739574844241159902101217660 & 2.022506 \\
30 & 6965244882542454937020619818702059053741377290255068284 & 2.024269 \\
31 & 521789556367416753405244328934612259884552790211207421244 & 2.025929 \\
32 & 39111402471791798530057405675011481922023421912454706904712 & 2.027493 \\
\multicolumn{1}{c}{\vdots} & \multicolumn{1}{c}{\vdots} &
\multicolumn{1}{c}{\vdots} \\
\\
60 & 
\parbox{.7\textwidth}{%
140299226506605674595297010256056870425972772255987818065754\\
8248104578242160531202625665330141151767257855947224
}\vspace{.1cm}
& 2.051236 \\

      \hline
    \end{tabular}
  }
  \caption{Results for $m=3$ (up to $n=10$ 
     by Aichholzer~\cite{Aichholzer-WWW}.)}
  \label{tab:m3}
\end{table}

\begin{table}[ht]
  \centering
  {\footnotesize
    \begin{sideways}
      \begin{tabular}{rrr}
        \hline
        \multicolumn{1}{c}{$n$\hspace{10mm}} &
        \multicolumn{1}{c}{$\#$ unimodular triangulations of $P_{4,n}$} &
        \multicolumn{1}{c}{capacity} \\
        \hline
        4 & 736983568 & 1.841066 \\
5 & 208902766788 & 1.880202 \\
6 & 61756221742966 & 1.908818 \\
7 & 18792896208387012 & 1.930751 \\
8 & 5831528022482629710 & 1.948080 \\
9 & 1835933384812941453312 & 1.962138 \\
10 & 584455230176565718869688 & 1.973785 \\
11 & 187686028049755013528577884 & 1.983601 \\
12 & 60685901262618326775192700244 & 1.991986 \\
13 & 19731268926382148037209063600412 & 1.999235 \\
14 & 6444884828545542240332780129017164 & 2.005567 \\
15 & 2113222804656668311309302902100087020 & 2.011147 \\
16 & 695163898467233943317499868644974218294 & 2.016103 \\
17 & 229316915701559537858641762255000442720116 & 2.020535 \\
18 & 75827610389461537709077484409103543784585710 & 2.024522 \\
19 & 25126215170054967918515172517611569017605311400 & 2.028130 \\
20 & 8341120564526486621411516194118239406742548614820 & 2.031409 \\
21 & 2773492141234111587868660829757660194821668552146720 & 2.034405 \\
22 & 923542836393484391634763396019705303948974557455679944 & 2.037151 \\
23 & 307927840785137540445620416718491375007439099052396877524 & 2.039680 \\
24 & 102788976776576952654837601573736966649525369087792041774754 & 2.042015 \\
25 & 34348007173983906887977536976143285540616055151406837464386952 & 2.044178 \\
26 & 11488735157368185493550332714289865763761071970295733672739011182 & 2.046188 \\
27 & 3846114798702837748429637799631405899720860204464245538618538080224 & 2.048061 \\
28 & 1288605098545745425026895981169173452565572617561190100930952999826626 & 2.049811 \\
29 & 432053057804721606144370507673814019049397228524243323674624721189513152 & 2.051449 \\
30 & 144960362040664009850310097114550102266012390064271163400409174321129109970 & 2.052986 \\
31 & 48666988587816100417429461522083834520238266167053837281579839713435133670880 & 2.054431 \\
32 & 16348321592766160525928861562545873132577437301251385744803248045014265044217096 & 2.055792 \\

        \hline
      \end{tabular}
    \end{sideways}
  }
  \caption{Results for $m=4$ (up to $n=8$ 
      by Aichholzer~\cite{Aichholzer-WWW}.)}
  \label{tab:m4}
\end{table}

\begin{table}[ht]
  \centering
  {\footnotesize
  \begin{tabular}{rrr}
    \hline
    \multicolumn{1}{c}{$n$\hspace{10mm}} &
    \multicolumn{1}{c}{$\#$ unimodular triangulations of $P_{5,n}$} &
    \multicolumn{1}{c}{capacity} \\
    \hline
    1 & 252 & 1.595455 \\
2 & 182132 & 1.747462 \\
3 & 182881520 & 1.829755 \\
4 & 208902766788 & 1.880202 \\
5 & 260420548144996 & 1.915513 \\
6 & 341816489625522032 & 1.941533 \\
7 & 464476385680935656240 & 1.961547 \\
8 & 645855159466371391947660 & 1.977388 \\
9 & 913036902513499041820702784 & 1.990240 \\
10 & 1306520849733616781789190513820 & 2.000871 \\
11 & 1887591165891651253904039432371172 & 2.009821 \\
12 & 2747848427721241461905176361078147168 & 2.017461 \\

    \hline
  \end{tabular}
  }
  \caption{Results for $m=5$ (up to $n=6$ 
    by Aichholzer~\cite{Aichholzer-WWW}.)}
  \label{tab:m5}
\end{table}

\begin{table}[ht]
  \centering
  {\footnotesize
  \begin{tabular}{rrr}
    \hline
    \multicolumn{1}{c}{$n$\hspace{10mm}} &
    \multicolumn{1}{c}{$\#$ unimodular triangulations of $P_{6,n}$} &
    \multicolumn{1}{c}{capacity} \\
    \hline
    1 & 924 & 1.641958 \\
2 & 2801708 & 1.784822 \\
3 & 12244184472 & 1.861743 \\
4 & 61756221742966 & 1.908818 \\
5 & 341816489625522032 & 1.941533 \\
6 & 1999206934751133055518 & 1.965553 \\
7 & 12169409954141988707186052 & 1.984082 \\

    \hline
  \end{tabular}
  }
  \caption{Results for $m=6$.}
  \label{tab:m6}
\end{table}


\end{document}